\documentclass[12pt]{amsart}
\title[The second order casimirs]
{Second order casimirs for the affine Krichever--Novikov algebras
$\widehat{\mathfrak{gl}}_{g,2}$ and $\widehat{\mathfrak{sl}}_{g,2}$}

\author[O. Sheinman]{O. K. Sheinman}

\address{Department of Geometry and Topology, Steklov Mathematical Institute,
Gubkina ul.~8, 117966, Moscow GSP-1, Russia\quad\textup{and}\endgraf
Independent University of Moscow}
\email{oleg@sheinman.mccme.ru}

\thanks{Supported in part by the RFBR Grant 99-01-00198}

\makeatletter
\def\subsection{\@startsection{subsection}{2}%
  \z@{.5\linespacing\@plus.7\linespacing}{-.5em}%
  {\def\@secnumpunct{\enspace}\normalfont\bfseries}}
\makeatother

\newcommand{\g}{\mathfrak{g}}

\let\a\alpha
\let\l\lambda
\let\ga\gamma
\let\w\omega
\let\ep\epsilon
\let\b\beta

\newcommand{\wh}[1]{\widehat{\vrule width0pt height1.4ex\smash[t]{\mathfrak{#1}}}}

\newcommand{\Z}{\mathbb Z}

\newcommand{\A}{\mathcal A}
\newcommand{\D}{\mathcal D^{1}}
\renewcommand{\L}{\mathcal L}
\newcommand{\M}{\mathcal M}
\newcommand{\U}{\mathcal U}

\newcommand{\eh}{\hat e}
\newcommand{\gh}{\hat\g}
\newcommand{\Gh}{\hat\g}
\newcommand{\Lh}{\widehat\L}

\newcommand{\At}{\widetilde\A}
\newcommand{\Lt}{\widetilde\L}

\newcommand{\GL}{\mathrm{GL}}
\newcommand{\tr}{\operatorname{tr}}
\newcommand{\ord}{\operatorname{ord}}
\newcommand{\res}{\operatorname{res}}
\newcommand{\rank}{\operatorname{rank}}

\newcommand{\rpi}[1]{\pi_{#1}}
\newcommand{\ainf}{\mathfrak{a}_\infty}
\newcommand{\vac}{|0\rangle}
\newcommand{\nord}[1]{{:}#1{:}}
\newcommand{\ldot}{\,.\,}

\DeclareSymbolFont{scripts}{U}{eur}{m}{n}
\DeclareSymbolFontAlphabet{\mathscr}{scripts}
\renewcommand{\k}{\mathscr{k}}
\newcommand{\ce}{\mathscr{c}}

\numberwithin{equation}{section}

\newtheorem{proposition}{Proposition}[section]
\newtheorem{conjecture}{Conjecture}[section]
\newtheorem{lemma}{Lemma}[section]
\newtheorem{theorem}{Theorem}[section]
\newtheorem{corollary}{Corollary}[section]

\theoremstyle{definition}
\newtheorem{example}{Example}[section]
\newtheorem{definition}{Definition}[section]

\theoremstyle{remark}
\newtheorem{remark}{Remark}[section]

\begin{document}
\begin{abstract}
The second order casimirs for the affine Krichever--Novikov algebras
$\widehat{\mathfrak{gl}}_{g,2}$ and $\widehat{\mathfrak{sl}}_{g,2}$
are described. More general operators which we call semi-casimirs are
introduced. It is proven that the semi-casimirs induce well-defined
operators on conformal blocks and, for a certain moduli space of
Riemann surfaces with two marked points and fixed jets of local
coordinates, there is a natural projection of its tangent space onto
the space of these operators. It is (non-formally) explained how
semi-casimirs appear in course of operator quantization of the second
order Hitchin integrals.
\end{abstract}
\maketitle

% \define\kazint{1}   Introduction
% \define\kazaff{2}   Affine Krichever-Novikov alg and their repr
% \define\kazvect{3}  Krichever-Novikov algebras of vector fields and their repr
% \define\kazmod{4}   Casimirs and the mod space ${\Cal M}_{g,2}$
% \define\kazapx1{5}

% 1
\section{Introduction}

The description of Casimir operators (casimirs, laplacians) is one of the
central questions of the representation theory of any Lie algebra. It is
difficult to list all the applications of casimirs.  The theory of special
functions as the eigenfunctions of casimirs, the construction of Hamiltonians
possessing symmetries and integrable systems, the investigation of spectral
properties of systems with symmetries are some of these applications. The second
order casimirs are of special interest in all these questions. In what follows
``casimir" always means ``second order casimir".

The \emph{casimirs} for a Lie algebra $\g$ can be characterized as the
operators which
\begin{enumerate}
\item commute with all the operators $\rho(g)$ for any (perhaps, from a given
class) representation $\rho$ of $\g$ and $g\in\g$,
\item can be expressed via $\rho(g)$ in a certain way.
\end{enumerate}
The casimirs belong to the wider class of the \emph{intertwinning operators},
which can be defined by omitting the second requirement (i.\,e., by imposing
only the requirement of commutativity).  For finite-dimensional semi-simple Lie
algebras, the investigation of casimirs is mainly based on I.\,M.\,Gelfand's
theorem about the center of the universal enveloping algebra. The developement
of the theory of Kac--Moody algebras suggested a new approach to constructing
casimirs \cite{KaRa}. This idea is closely related to the following fundamental
phenomenon.  With each (so-called \emph{admissible}) representation of an affine
Kac--Moody algebra $\gh$ a certain representation of the Virasoro algebra
($\mathit{Vir}$) is canonically associated. The latter is called the Sugawara
representation. It acts in the same linear space as the initial representation
of $\gh$.  Denote by $\D$ the sum of $\mathit{Vir}$ and $\gh$ with their centers
identified.  In an admissible representation of $\D$, each element $e$ of
$\mathit{Vir}$ acts first in the usual way and second via its Sugawara operator.
For a semi-simple $\g$ and its loop algebra $\gh$, a certain linear combination
$\Delta_e$ of these two actions always produces an intertwinning operator. One
of them is called the \emph{casimir}, namely the one that corresponds to the
vector field of zero degree (here we do not consider the degenerate case of
so-called critical level, in which there is an infinite number of casimirs of
this form).

Krichever and Novikov introduced and investigated \cite{KNFa,KNFb,KNFc} a
certain type of tensors on Riemann surfaces, namely the meromorphic tensors that
have poles only at two fixed points $P_\pm$ of the Riemann surface.  We call
such tensors Krichever--Novikov tensors. For example, we speak of
Krichever--Novikov functions (which form the associative algebra $\A$),
Krichever--Novikov vector-fields, $\l$-forms, etc.  Two new classes of Lie
algebras were introduced in \cite{KNFa,KNFb,KNFc}. These are the centrally
extended Lie algebras of Krichever--Novikov vector-fields and $\g$-valued
Krichever--Novikov functions, where $\g$ is a finite-dimensional complex
semisimple or reductive Lie algebra.  Let $\smash{\Lh}$ denote the first of them (it is
called the \emph{Virasoro-type} Krichever--Novikov algebra) and $\gh$ denote the
second one (it is called the \emph{affine type} Krichever--Novikov algebra). For
zero genus, they coincide with the usual Virasoro algebra and affine Kac--Moody
algebra, respectively. For affine Krichever--Novikov algebras, the
above-mentioned idea of constructing casimirs as the $\Delta_e$-operators was
first realized in \cite{SSS}.

This paper is devoted to the description of the second order casimirs (and a
certain generalization) for affine Krichever--Novikov algebras. We only consider
Krichever--Novikov algebras which correspond to $\g=\mathfrak{sl}(l)$ (the
semi-simple case: properties of simplicity are of no importance here, see
Lemma~4.1 below) and $\g=\mathfrak{gl}(l)$ (the reductive case).  We use the
more detailed notation for $\gh$ in these cases: $\smash{\wh{sl}_{g,2}}$
for the first case and $\smash{\wh{gl}_{g,2}}$ for the second one. Here
$g$ stands for the genus of the Riemann surface in question and the index~$2$
corresponds to the number of marked points.

Our description of casimirs is based on the above construction of the
$\Delta_e$'s where $\mathit{Vir}$ and the affine Kac--Moody algebra are replaced
by their Krichever--Novikov counterparts $\smash{\Lh}$ and $\gh$.  The questions that
immediately arise are as follows: How many independent (second order) casimirs
exist? Why only one of the intertwinning operators described above is considered
as a casimir for the Kac--Moody algebra? The usual argumentation is that one of
the elements of $\mathit{Vir}$ is distinguished because it defines a grading of
the affine algebra. This doesn't work for Krichever--Novikov case, because there
is no reason to distinguish any element of $\L$ here (in particular,
Krichever--Novikov algebras are not graded). In what follows we relate the
question about the number of independent casimirs with a certain cocycle $\ga$
on~$\D$. It turns out that only the elements $e\in\L$ such that $\ga(e,A)=0$ for
any $A\in\A$ produce casimirs.  Using this, we prove that \emph{the affine
Krichever--Novikov algebras $\widehat{\mathfrak{sl}}_{g,2}$ and
$\widehat{\mathfrak{gl}}_{g,2}$ possess only one second order casimir}. In
particular, this holds for Kac--Moody algebras ($g=0$).

Further, we also consider weaker conditions than those which define casimirs.
Let $\A_\pm\subset\A$ be the subalgebra consisting of elements of positive
order at the point $P_\pm$. For a certain subspace $\A_0\subset \A$ (see the
definition below) the following decomposition holds \cite{KNFa}:
$\A=\A_-\oplus\A_0\oplus\A_+$.
%Define the subspace $\At_-\subset\A$ by the formula
%$\At_-=\A _-\oplus\A _0$.
Consider vector fields $e$  such that $\ga(e,A)=0$ for any $A\in\A_-\oplus\A_0$.
We call the operators $\Delta_e$ for such vector fields \emph{semi-casimirs}.

It turns out to be that there is an interesting geometrical connection between
semi-casimirs and certain moduli spaces of Riemann surfaces. Let
$\M_{g,2}^{(p)}$ be the moduli space of Riemann surfaces of genus $g$ with
$2$ marked points $P_\pm$ and fixed jets of local coordinates (of order $1$ for
$P_+$ and of order $p$ for $P_-$).  For a point $\Sigma\in\M_{g,2}^{(p)}$, let
$T_\Sigma\M_{g,2}^{(p)}$ denote the tangent space to this moduli space
at~$\Sigma$.  Consider also the space of co-invariants (or \emph{conformal
blocks}) of the \emph{regular subalgebra}~$\g_r$ (see Section~4(c) for the
definition) in some representation of $\gh$ for $\g=\mathfrak{gl}(l)$.  It turns
out to be that semi-casimirs are well defined on  conformal blocks and (by
Theorem~4.2 below) \emph{for a proper~$p$ there is a natural projection of
$T_\Sigma\M_{g,2}^{(p)}$ onto the space of operators induced by semi-casimirs on
conformal blocks over $\Sigma$}.

The paper is organized as follows.  In Section~2 we systematically introduce
affine Krichever--Novikov algebras and their representations. For the proofs we
refer to~\cite{Shf}. In Section~3 we introduce the Krichever--Novikov algebras
of vector fields. We also define the Sugawara representation and consider its
commutation relations with the the representation of the affine algebra. These
two sections can be regarded as an introduction to Krichever--Novikov algebras
and their representations. In Section~4 we introduce casimirs and semi-casimirs
and formulate the results mentioned above (Theorems~4.1 and~4.2). In Section~5
we restate some of the results of Section 4\ from another point of view.  Note
that in Section~4 we do not use the classification theorem for 2-cocycles on
$\D$. Such a theorem is proved in \cite{ArCoKaPr} for genus zero.  As for
positive genus, we don't know any published proof of the theorem. Nevertheless
the common opinion is that the result is true\footnote
{In a private communication, B.\,Feigin mentioned that he has a proof of
the result.}.
To demonstrate that both approaches give the same results, in Section~5 we prove
some results of Section~4 assuming the classification of 2-cocycles on $\D$ to
be known.

The setting of the problems and some approaches in this paper were significantly
stimulated by my joint work with M.\,Schlichenmaier. In particular for a
commutative $\g$ he realized that the cocycles on $\D$ are actually obstructions
for the $\Delta_e$'s to be casimirs (see Lemma~4.2).  I am thankful to
M.\,Schlichenmaier for stimulating discussions and for the hospitality at the
University of Mannheim. I am also thankful to B.\,Feigin and S.\,Loktev for
discussions about cocycles on $\D$ and co-invariants.

\clearpage

% 2
\section{Affine Krichever--Novikov algebras and their representations}

% (a)
\subsection{Affine Krichever--Novikov algebras}

Let $\Sigma$ be a compact algebraic curve over $\mathbb{C}$ of genus $g$ with
two marked points $P_{\pm}$, $\A (\Sigma,P_\pm)$ be the algebra of meromorphic
functions on $\Sigma$ which are regular outside the points $P_\pm$, ${\g}$ be a
complex semi-simple or reductive Lie algebra.  Then
\begin{equation}
\gh=\g\otimes_{\mathbb{C}}\A (\Sigma,P_\pm)\oplus\mathbb{C}c
                                             \label{alg}
\end{equation}
is called the \emph{affine Krichever--Novikov algebra}
\cite{KNFa, Shea}. The bracket on $\gh$ is given by the relations
\[
[x\otimes A,y\otimes B]=[x,y]\otimes AB+\ga(x\otimes A,y\otimes B)c , \qquad
      [x\otimes A,c]=0,
\]
where $\ga$ is the cocycle defined by the formula
\begin{equation}
\gamma (x\otimes A,y\otimes B)=(x,y)\res_{P_+}(A\,dB);
                                                        \label{cocycle}
\end{equation}
here $(\mspace{2mu}\cdot\,{,}\cdot\mspace{2mu})$ is a nondegenerate invariant
bilinear form on $\g$ (e.\,g.\
for $\g=\mathfrak{gl}(l)$ we take $(x,y)=\tr (xy)$). We will usually suppress
the symbol $\otimes$ in our notation.  We also often write $\A$ instead of
$\A(\Sigma,P_\pm)$.

In \cite{Shf}\ the fermion (or wedge) representations of affine
Krichever--Novikov algebras were introduced. They form a rather representative
class.  Here we use them as a basic example and a model for our
constructions. In the remainder of this section, we systematically describe these
representations but refer to \cite{Shf} for proofs.

% 2(b)
\subsection{Holomorphic bundles. Tjurin parameters}

Each fermion representation is related to a holomorphic vector bundle on the
Riemann surface $\Sigma$. For this reason consider a holomorphic bundle $F$ of
rank $r$ and degree $gr$ on $\Sigma$.  By the Riemann--Roch theorem, the
bundle~$F$ possesses $r$ holomorphic sections $\Psi_1,\dots,\Psi_r$ which form a
basis over each point except for $gr$ points.  For a generic situation (which is
only considered here), one can choose these points to be mutually different. We
call these points the \emph{degeneration points} and denote them
$\ga_1,\dots,\ga_{gr}$.  Following the terminology of \cite{KNFd,KNU} we call
the set of sections introduced a \emph{framing} and we call a bundle with a
given framing a \emph{framed bundle}.

In a local trivialization of $F$ one can consider the sections
$\Psi_1,\dots,\Psi_r$ as vector-valued functions (say columns).  These functions
can be arranged into a matrix $\Psi$.  This matrix is invertible everywhere
except for the points $\ga_i$, $i=1,\dots,gr$, in which $\det\Psi$ has simple
zeroes: $\det\Psi(\ga_i)=0$, $(\det\Psi)'(\ga_i)\ne 0$.  We call
$D=\ga_1+\dots+\ga_{gr}$ the \emph{Tjurin divisor} of the bundle~$F$.  Impose
one more condition of genericity: $\rank \Psi(\ga_i)=r-1$, $i=1,\dots,gr$.  Then
for each $i=1,\dots,gr$ a non-trivial solution to the system of linear equations
$\Psi(\ga_i)\a_i=0$ exists and is unique up to a scalar factor.  Let us
introduce notation $\a_i=(\a_{ij})^{j=1,\dots,r}$ ($i=1,\dots,gr$).  The Tjurin
divisor of the bundle and the numbers $(\a_{ij})^{j=1,\dots,r}_{i=1,\dots,gr}$
are called the Tjurin parameters of the bundle $F$ \cite{KNFd,KNU}.  The framing
is defined uniquely up to action of the group $\GL(r)$ on $\Psi$ by
\emph{right} multiplication in contrary with gluing functions which act
\emph{on the left}.  We see that this $\GL(r)$-action commutes with the action
of gluing functions; hence it sends sections into sections. This is why the
Tjurin parameters are determined uniquely up to a scalar factor and \emph{left}
$\GL(r)$-action. Let us emphasize that equivalent framed bundles have the same
set $\ga_1,\dots,\ga_{gr}$ while for non-framed bundles only the class of the
divisor is invariant. According to \cite{Tur} the Tjurin parameters determine the
bundle uniquely up to equivalence.

In \cite{KNFd,KNU} the following description of the space of meromorphic
sections of the bundle $F$ in terms of its Tjurin parameters is proposed (only
those meromorphic sections are considered which are holomorphic outside the
points~$P_\pm$).  In the fibre over an arbitrary point $P$ outside the support
of the divisor $D$ the elements $\Psi_j(P)$ ($j=1,\dots,r$) form a base.  Hence
for each meromorphic section $S$ its value $S(P)$ can be expanded in terms of
this basis.  Thus, to each section $S$ one can assign a vector-valued function
$f=(f_1,\dots,f_r)^T$ on the Riemann surface $\Sigma$ so that outside $D$ one has
\begin{equation}
S(P)=\sum\limits_{j=1}^r\Psi_j(P)f_j(P).      \label{razl}
\end{equation}
By Cramer's rule
$f_j=\det(\Psi_1,\dots,\Psi_{j-1},S,\Psi_{j+1},\dots,\Psi_r)\,
(\det\Psi)^{-1}$.  This shows that the functions $f_j$ can be continued to the
points of the divisor~$D$. They will have there at most simple poles because
$\Psi_1,\dots,\Psi_r$ are holomorphic at the points of the divisor $D$ and
$\det\Psi$ has simple zeroes there. By \eqref{razl} in local coordinates in a
neighborhood of a point $\ga_i$ one has $S(z)=\Psi(\ga_i)(\res_{\ga_i}f)z^{-1}+
O(1)$.  The left-hand side of the latter relation is holomorphic. Hence the
residues of the functions~$f_j$, $j=\nobreak 1,\dots,r$, at the point $\ga_i$ satisfy the
system of linear equations $\Psi(\ga_i)(\res_{\ga_i}f)=\nobreak0$. This is just the
system which determines the Tjurin parameters at the point $\ga_i$.  By
assumption the rank of the matrix $\Psi(\ga_i)$ is equal to $r-1$. Hence the
vectors $\res_{\ga_i}f$ and $\a_i$ are proportional.

% Proposition 2.1
\begin{proposition}[\cite{KNFd,KNU}]
Consider the space of meromorphic sections of the bundle $F$ which are
holomorphic outside the marked points $P_\pm$. This space is isomorphic to the
space of meromorphic vector-valued functions $f=(f_1,\dots,f_r)^T$ \upn(on the same
Riemann surface\upn) which are holomorphic outside the points $P_\pm$ and the
divisor~$D$, have at most simple poles at the points of the $D$ and satisfy the
relations
\[
(\res_{\ga_i}f_j)\mspace{2mu}\a_{ik}=(\res_{\ga_{i}}f_k)\mspace{2mu}\a_{ij},\qquad
      i=1,\dots,gr,\ j=1,\dots,r.
\]
\end{proposition}

Let us denote the space of vector-valued functions just defined by $F_{KN}$.

% 2(c)
\subsection{The Krichever--Novikov bases}

Let us introduce a basis in $F_{KN}$ having in mind constructing semi-infinite
monomials on this space. For each pair of integers $n,j$ ($0\le j<r$) let us
construct a vector-valued function $\psi_{n,j}\in F_{KN}$. The integer $n$ is
called the \emph{degree} of~$\psi_{n,j}$. The function $\psi_{n,j}$ is specified
by its asymptotic behavior at the points $P_{\pm}$. Consider $\psi_{n,j}$ as a
column and arrange the matrix $\Psi_n$ of these columns.  We require
\begin{equation}
\Psi_n(z_+)
   =z_+^n \sum\limits_{s=0}^\infty \xi^+_{n,s}z_+^s,\qquad
    \xi^+_{n,0}=\begin{pmatrix}
       1   &    *   & \dots &    *   \\
       0   &    1   & \dots &    *   \\
 \hdotsfor{4} \\
       0   &    0   & \dots &    1
    \end{pmatrix},                                  \label{bas1}
\end{equation}
and
\begin{equation}
\Psi_n(z_-)
   =z_-^{-n}\sum\limits_{s=0}^\infty \xi^-_{n,s}z_-^s,\qquad
    \xi^-_{n,0}=\begin{pmatrix}
       *   &    0   & \dots &    0   \\
       *   &    *   & \dots &    0   \\
\hdotsfor{4}\\
       *   &    *   & \dots &    *
    \end{pmatrix},                                  \label{bas2}
\end{equation}
where $z_\pm$ is a local coordinate at the point $P_\pm$ and $*$'s denote
arbitrary complex numbers.

Thus the matrix $\Psi_n$ has a zero of order $n$ at one of the points $P_\pm$
and the pole of order $n$ at the other one. Its determinant has $gr$ simple
poles at the points of the divisor $D$ and an additional apriory not fixed
divisor of zeroes outside the points $P_\pm$.  We call $\{\psi_{n,j}\}$ the
\emph{Krichever--Novikov basis} in $F_{KN}$.

The case $r=1$ is exceptional and needs special prescriptions for the
Krichever--Novikov bases.

% Example 2.1
\begin{example}
The prescription for the Krichever--Novikov basis in the algebra $\A
(\Sigma,P_\pm)$ (that is for a rank 1 trivial bundle $F$) as introduced in
\cite{KNFa} is as follows:
\begin{equation}
A_m=\a_m^{\pm}z_{\pm}^{\pm m+\varepsilon_\pm}(1+O(z_\pm)),\qquad
     \a_m^\pm\in\mathbb{C}, \ \a_m^+=1,         \label{class}
\end{equation}
where $\varepsilon_+=0$ for any $m\in\mathbb{Z}$, $\varepsilon_-=-g$ for $m>0$
or $m<-g$, and $\varepsilon_-=-g-1$ for $-g\le m\le 0$. For $m>0$ or $m<-g$ the
sum of orders at the marked points is equal to $(-g)$ (notice that $r=1$
in this case).  Therefore there are exactly $g$ zeroes (and no poles) outside
$P_\pm$.

Denote by $\A_+$ (respectively, $\A_-$, $\A_0$) the linear space spanned
by $A_m$'s, $m>0$ (respectively, $m<-g$, $-g\le m\le 0$). The
following decomposition holds \cite{KNFa}:
\begin{equation}
\A=\A_-\oplus\A_0\oplus\A_+ .                   \label{decomp}
\end{equation}
\end{example}

Another example is given in Section~3(a).

% Proposition 2.2
\begin{proposition}[\cite{Shf}]
$1^\circ.$ There exists a unique matrix-valued function $\Psi_n$ which satisfies
conditions \eqref{bas1}, \eqref{bas2}.

$2^\circ.$ The dimension of the space generated by the vector-valued functions
$\psi_{n,j}$ \upn($n$ being fixed\upn) is equal to $r$.
\end{proposition}

The algebra $\A(\Sigma,P_\pm)$ naturally acts in the space $F_{KN}$ multiplying
its elements by functions. We introduce the structure constants of that action
by means of the following relation \cite[Proposition~2.3]{Shf}:
\begin{equation}
A_m \psi_{n,j} =\sum\limits_{k=m+n}^{m+n+{\bar g}}
                     \sum\limits_{j'=0}^{r-1}
                     C^{k,j'}_{m,n,j} \psi_{k,j'},
                                               \label{act}
\end{equation}
where $\bar g=g+1$ if $-g\le m\le 0$ and ${\bar g}=g$ otherwise.  This relation
expresses the fact that the $\A(\Sigma,P_\pm)$-module $F_{KN}$ is \emph{almost
graded}, i.\,e., $k$ in \eqref{act} satisfies $|m+n-k|<\mathit{const}$ and the
constant does not depend on $m$,~$n$.

Let $\tau$ be the representation of the algebra $\g$ in the linear space
$\mathbb{C}^l$. We write $\tau(x)=(x^i_{i'})$, where $x\in\g$, $(x^i_{i'})\in
\mathfrak{gl}(l)$ is the matrix which represents the element $x$, and the indices
$i$ and~$i'$ run over $\{1,\dots,l\}$. To each basis element $\psi_{n,j}$
($n\in\mathbb{Z}$, $j=0,\dots,r-1$) let us assign a set of basis elements
$\{\psi_{n,j}^i\colon i=1,\dots,l\}$ of $F_{KN}^\tau=F_{KN}\otimes\mathbb{C}^l$.
Define an action of the Lie algebra $\g\otimes\A(\Sigma,P_\pm)$ on $F_{KN}^\tau$:
\begin{equation}
(xA_m)\psi_{n,j}^i = \sum\limits_{k=m+n}^{m+n+{\bar g}}
\sum\limits_{j'=0}^{r-1}\sum\limits_{i'}
               C_{m,n,j}^{k,j'}x^i_{i'}\psi_{k,j'}^{i'}.
                                              \label{actg}
\end{equation}
The action of $\g$ does not change the indices $n$, $j$; and, for given $n$ and
$j$, it is $\psi^l_{n,j}$ that plays the role of the highest weight vector. It
follows from the definitions that this action is almost graded.

% 2(d)
\subsection{Fermion representations with highest weight}

Let us enumerate the symbols $\psi_{n,j}^i$ in the lexicographic order of the
triples $(n,j,i)$. Let $N=N(n,j,i)$ be the number of a triple $(n,j,i)$. We
normalize this unmbering by the condition $N(-1,0,l)=0$. Introduce the following
notation: $\psi_N=\psi_{n,j}^i$.

Introduce the fermion representation which corresponds to $F$, $\tau$, as
follows.  The space $V_F$ of the representation is generated over $\mathbb C$ by
the formal expressions (\emph{semi-infinite monomials}) of the form
$\psi_{N_0}\wedge\psi_{N_1} \wedge\dots$, where $N_0<N_1<\dots$ and the sign of
the monomial changes under the transposition of $\psi_N$, $\psi_{N'}$.  We also
require that $N_k=k+m$ for $k$ sufficiently large ($k\gg 1$). Following
\cite{KaRa} we call $m$ the \emph{charge} of the monomial.

For a monomial $\psi$ of charge $m$ the degree of $\psi$ is defined as follows:
\begin{equation}
\deg\psi =\sum\limits_{k=0}^\infty (N_k-k-m).  \label{mdeg}
\end{equation}
Observe that there is an arbitrariness in the numbering of the $\psi_{n,j}^i$'s
for $n$ fixed; the degree of a monomial just defined does not depend on this
arbitrariness.

The representation $\rpi{F,\tau}$ in the linear space $V_F$ is defined as
follows. By \eqref{actg} each element of $\gh$ acts on the symbols $\psi_N$ by
linear substitutions in an almost graded way. Moreover, the number of the
symbols $\psi_N$ of a fixed degree does not depend on this degree. This exactly
means that if the symbols $\psi_N$ are considered as a formal basis of an
infinite-dimensional linear space then the action of an element of the algebra
$\gh$ can be given by an infinite matrix with only a finite number of nonzero
diagonals in this basis. Following \cite{Kac}, we denote the algebra of such
matrices by~$\ainf$.

% Remark 2.1
\begin{remark}
In other terms, $\ainf$ is the algebra of difference operators on a
$1$-di\-men\-sional lattice. This remark yields a relation to the results of
\cite{KNRm}.
\end{remark}

Thus, fixing the basis $\{\psi_N\}$ provides us with an imbedding of
$\g\otimes\A$ into $\ainf$. Since $\ainf$ has a standard action in $V_F$, we
obtain a representation of $\g\otimes\A$ in $V_F$.  Recall \cite{KaRa, Shf} that
the action of the basis element $E_{IJ}\in\ainf$ on a semi-infinite monomial
$\psi=\psi_{I_0}\wedge\psi_{I_1}\wedge\dotsb$ is defined by the Leibniz rule:
\begin{equation}
 r(E_{IJ})\psi=(E_{IJ}\psi_{I_0})\wedge\psi_{I_1}\wedge\dots+
                   \psi_{I_0}\wedge(E_{IJ}\psi_{I_1})\wedge\dots+
                   \dotsb .
                                            \label{infact}
\end{equation}
Due to the condition $I_k=k+m$ ($k\gg 1$) the action \eqref{infact} is
well defined for $I\ne J$. For $I=J$, the right-hand side of \eqref{infact}
contains infinitely many terms. In this case the standard regularization is used
\cite{KaRa, Shf}; it results in a projective representation $\rpi{F,\tau}$ of
$\g\otimes\A$.  This is an almost graded representation, as it follows from
almost-gradedness of the action \eqref{actg} in the space $F_{KN}$.  Another
standard procedure enables us to transform $\ga$ to any cohomological cocycle.
This can be achieved by adding certain scalars to the operators
$\rpi{F,\tau}(x_\a A_m)$ where $x_\a$ are the generators of~$\g$.  Such
modification of the representation results in adding a coboundary to
$\ga$. Hence, by \cite{PS} $\ga$ can be written as $(2\pi i)^{-1}(x,y)\oint
A\,dB$ where the notation is the same as in~\eqref{cocycle}. A 2-cocycle on
$\g\otimes\A$ is called \emph{local} if there exist $L\in\mathbb{Z}_+$ such that
$\ga(xA_i,yA_j)=0$ for any $i,j\in\mathbb{Z}$ such that $|i-j|>L$ and any
$x,y\in\g$.  Since $\ga$ is obviously local, by \cite{KNJGP,KNfl} it is
cohomological to the multiple of the cocycle \eqref{cocycle}. Therefore we can
consider $\rpi{F,\tau}$ as a representation of $\gh$ without loss of generality.

% Conjecture 2.1 ?
\begin{conjecture}
The equivalence classes of fermion representations of $\gh$ are in one-to-one
correspondence with the following pairs\upn: an equivalence class of a
holomorphic bundle of rank~$r$ and degree $gr$ on $\Sigma$ and an equivalence
class of an $l$-dimensional representation of the algebra~$\g$.
\end{conjecture}

For $l=r$ this conjecture is proven in \cite{Shf} (for the definition of equivalence of
fermion representations see also \cite{Shf}).

Obviously, the action of $\ainf$ respects the charge of a monomial because the
infinite ``tails" of the monomials both on the right-hand and on the left-hand
sides of \eqref{infact} are the same.  Hence the fermion space can be decomposed
into a direct sum of $\gh$-invariant subspaces of certain charges. For
$\g=\mathfrak{gl}(l)$, the space of charge $m$ is generated by the monomial
$\vac=\psi_m\wedge\psi_{m+1}\wedge\dotsb$ under the action of the universal
enveloping algebra $\U(\gh)$.  This monomial is called a \emph{vacuum}
(\emph{monomial}) of charge~$m$.  Each vacuum monomial possesses the following
property: $\rpi{F,\tau}(xA)\,\vac=0$ for $A\in\A_+$ and also for $A=1$ and any
strictly upper-triangular matrix $x\in\g$.

% 3
\section{Krichever--Novikov algebras of vector fields
and~their~representations}

% 3(a)
\subsection{Algebras of vector fields and their central extensions}

Let $\L$ be the Lie algebra of meromorphic vector fields on $\Sigma$ which are
allowed to have poles only at the points $P_\pm$ \cite{KNFa,KNFb,KNFc}.  As it
was first noticed in \cite{KNFa}, for each integer $s\ge 0$ there is a pair of
subalgebras $\L^{(s)}_\pm$ in $\L$ which consist of vector fields of order not
less than $s$ at the points $P_\pm$ respectively. The following decomposition
into a direct sum of subspaces holds \cite{KNFa}:
$\L=\L^{(s)}_+\oplus\L^{(s)}_0\oplus\L^{(s)}_-$, where $\L^{(s)}_0$ is a
complementary space. In what follows we are mainly interested in the case $s=2$
in connection with deformations of moduli of Riemann surfaces.

Considered as a linear space, $\L$ has the Krichever--Novikov basis $\{e_m\colon
m\in\mathbb{Z}\}$. In the case $g>1$ the basis vector field $e_m$ is given by its
asymptotic behavior in the neighborhoods of the points $P_\pm$ \cite{KNFa}:
\begin{equation}
e_m(z_\pm)=\epsilon_m^\pm z_\pm^{\pm m+\varepsilon_\pm}
(1+ O(z_\pm))\frac{\partial}{\partial z_\pm},\qquad
    \epsilon_m^\pm\in\mathbb{C},\ \epsilon_m^+=1, \label{vect}
\end{equation}
where $\varepsilon_+=1$, $\varepsilon_-=1-3g$, $z_\pm$ is a local coordinate at
$P_\pm$. Thus the subalgebra $\L^{(s)}_+$ is generated by the elements $e_m$,
$m\ge s-1$, the subalgebra $\L^{(s)}_-$ by the elements $e_m$, $m\le -s-3g+1$
and the subspace $\L^{(s)}_0$ can be chosen to be generated by the
elements~$e_m$, $-s-3g+1<m<s-1$.

We will consider central extensions of the algebra $\L$. Each central extension
is given by a 2-cocycle of the form
\begin{equation}
\chi(e,f)=\res_{P_+}\biggl(\frac{1}{2}(e'''f-ef''')-R(e'f-ef')\biggr)
                                     \label{vcoc}
\end{equation}
(cf.\ Lemma~5.1) where $R$ is a projective connection, i.\,e., such a function
of a point and a local coordinate that transforms under a change of coordinates
as follows:
\[
R(u)u_z^2=R(z)+\frac{u_{zzz}}{u_z}-
    {\frac{3}{2}} \biggl(\frac{u_{zz}}{u_z}\biggr)^2.
\]
This ensures that the right hand side of \eqref{vcoc} is indeed the residue of a
well-defined 1-form. Let $\Lh$ denote a central extension of $\L$ given
by~\eqref{vcoc}.

% 3(b)
\subsection{The action of vector fields in the space $V_F$}

Consider the action of the algebra $\L$ in the space $F_{KN}$. Observe that we
cannot apply naively some vector field $e\in\L$ considered as a first order
differential operator to a vector-valued function $\psi\in F_{KN}$. The reason
is that a generic vector field $e\in\L$ is of order zero at the points
$\gamma_1,\dots,\gamma_{gl}$ of the Tjurin divisor of the bundle. Hence, $e\psi$
generically has poles of order two at these points. Thus, $e\psi\notin F_{KN}$.
Nevertheless an $\L$-action on $F_{KN}$ can be defined.

According to \cite{Bol} each holomorphic bundle, in particular $F$, can be
endowed with a meromorphic (therefore flat) connection $\nabla$ which has
logarithmic singularities at the points $P_\pm$. Since $\nabla$ is flat,
$[\nabla_e,\nabla_f]-\nabla_{[e,f]}=0$ for any $e,f\in\L$.  Hence
$\nabla_{[e,f]}= [\nabla_e,\nabla_f]$, and $\nabla$ defines a representation of
$\L$ in the space of holomorphic sections of $F$.  Being conjugated by the matrix
$\Psi$ (see Section~2(b)), this representation can be transferred to the space
$F_{KN}$ of Krichever--Novikov vector-valued functions. In what follows we fix an
arbitrary connection $\nabla$ which is logarothmic at $P_\pm$ and regular
elsewhere.  Applying the procedure described above (Section~2) of lifting
the representation from the space $F_{KN}$ to the space $V_F$ we obtain the
representation of the central extension of $\L$.

In local coordinates, let $e=E(z)\partial$ where
$\partial=\frac{\partial}{\partial z}$ and
$\nabla_e=E(z)(\partial+\omega)$. Then the action of the vector field $e$ on
$\psi$ can be written in the following two equivalent forms:
\begin{equation}
\begin{aligned}
e\psi &= \Psi^{-1}\nabla_e\Psi\psi,\\
e\psi &= E(\partial+\omega_\Psi)\psi,
\end{aligned}	\label{eact}
\end{equation}
where $\omega_\Psi=\Psi^{-1}\omega\Psi+\Psi^{-1}\Psi'$ ($\Psi$ being introduced
in Section~2(b)).  Due to \eqref{eact}, $e\mapsto\nobreak e\psi$ is a well-defined first
order differential operator in Krichever--Novikov vec\-tor-valued functions
(indeed, $(\partial+\omega_\Psi)\psi$ is a vector-valued 1-form; multiplication
by $E$ returns it into the space of vector-valued functions).

% Lemma 3.1
\begin{lemma}
$1^\circ$. The space $F_{KN}$ is invariant under the action \eqref{eact}.

$2^\circ$. The action \eqref{eact} is almost graded.
\end{lemma}

\begin{proof}
The $1^\circ$ follows from the definition of the action \eqref{eact}.

To prove $2^\circ$, observe that $\Psi$ is regular and nondegenerate at the
points $P_\pm$. Hence $\omega_\Psi$ also has simple poles at $P_\pm$. Hence
$\partial$ and $\omega_\Psi$ equally decrease the order of $\psi$ while acting
in local coordinates in neighborhoods of $P_\pm$.

Let us make use of \eqref{vect}.  For $n\ne 0$, the order of $e_m\psi_{n,j}$ is
equal to $n+m$ at $P_+$ and $-n-m-3g$ at $P_-$; for $n=0$, these orders are
equal to $n+m$ and $-n-m-3g$ respectively.  Hence the following relations hold
for some constants $D^{k,j'}_{m,n,j}$:
\begin{equation}
e_m \psi_{n,j}^i
     =\sum\limits_{k=m+n+\varepsilon}^{m+n+3g+\varepsilon}
                     \sum\limits_{j'=0}^{l-1}
                     D^{k,j'}_{m,n,j} \psi_{k,j'}^i,
                                               \label{vact}
\end{equation}
where $\varepsilon=0$ for $n\ne 0$ and $\varepsilon=1$ for $n=0$.
\end{proof}

Since the $\L$-action on $F_{KN}$ is almost graded, it can be continued
to the space of the fermion representation as representation of $\Lh$
(by the scheme described in detail in Section~2).

% Lemma 3.2
\begin{lemma}
The elements of the subalgebra $\L^{(2)}_+$ annihilate the vacuum
vectors of the fermion representations \upn(of highest weight\upn).
\end{lemma}

\begin{proof}
By definition (see Subsection (a)) $e_m\in\L^{(2)}_+$ if and only if $m\ge
1$. In this case, by \eqref{vact} the action of the vector field $e$ increases
the value of the index $n$: $k>n$ for all $\psi_{k,j'}^i$ which occur in
\eqref{vact} on the right-hand side. But since the vacuum monomial contains
$\psi_{n,j}^i$, it also contains \emph{all} the $\psi_{k,j'}^i$, $k>n$ on the
right of it. This is why $e_m$ annihilates the vacuum.
\end{proof}

% 3(c)
\subsection{Sugawara representation}

The following is one of the most fundamental facts of representation theory of
affine algebras. Each \emph{admissible} representation of an affine algebra
canonically defines a representation of the (corresponding) Virasoro type
algebra in the same space. The latter is called the Sugawara representation. For
Krichever--Novikov algebras the Sugawara representation is considered in
\cite{Bonora,KNFb,SSS}. Here we recall the basic facts about the Sugawara
representation.

Let us introduce the following notation. For $x\otimes A\in\Gh$, we denote the
corresponding operator of representation by $x(A)$. For a basis element
$A_n\in\A$, we denote $x(A_n)$ by $\ x(n)$.

A module (representation) $V$ over the Lie algebra $\gh$ is said to be
\emph{admissible} if $x(n)v=0$, for each $v\in V$, $x\in\g$ and $n$ sufficiently
large.

Let $V$ be an admissible module. It is assumed that the central element of $\gh$
acts by multiplication by a scalar $\ce\in\mathbb C$. In this paper $V$ can be thought
of as one of fermion representations which were introduced in Section~2.

From now on, till the end of this section, $\g$ is either simple or abelian Lie
algebra. The reason is that the existence of a nondegenerate invariant bilinear
form on $\g$ is crucial.  Any simple $\g$ has such a form; if $\g$ is abelian,
we require it to be equipped with such a form.

Take a basis $u_i$, $i=1,\dots,\dim\g$, of $\g$ and the corresponding dual basis
$u^i$, $i=1,\dots,\dim\g$, with respect to the invariant nondegenerate symmetric
bilinear form $(\mspace{2mu}\cdot\,{,}\cdot\mspace{2mu})$.  The Casimir element
$\Omega^0=\sum_{i=1}^{\dim \g} u_iu^i $ of the universal enveloping algebra
$U(\g)$ is independent of the choice of the basis. In what follows the summation
over $i$ is assumed in all cases when $i$ occurs in both super and subscripts.

Let $2\k$ be the eigenvalue of $\Omega^0$ in the adjoint representation. For a
simple $\g$, $\k$ is the \emph{dual Coxeter number}. In abelian case, $\k=0$.

Introduce bases in the spaces of Krichever--Novikov 1-forms and quadratic
differentials.  Denote their elements by $\{\omega^m|m\in\mathbb{Z}\}$ and
$\{\Omega^k\colon k\in\mathbb{Z}\}$ respectively. We define the basis elements
by the following \emph{duality relations}:
\[
\res_{P_+}A_m\omega^n=\delta^n_m,\qquad
\res_{P_+}e_m\Omega^k=\delta^k_m,
\]
where $\delta^n_m$ is the Kronecker symbol. The $\omega^n$'s have the following
asymptotic behavior \cite{KNFb}:
\begin{equation}
\w^n(z_\pm)=\mu_n^\pm z_\pm^{\mp n+\varepsilon_\pm}(1+O(z_\pm))dz_\pm,
\qquad \mu_n^\pm\in\mathbb{C},\ \mu_n^+=1, \label{diff}
\end{equation}
where $\varepsilon_+=-1$, $\varepsilon_-=0$, $z_\pm$ is a local
coordinate at $P_\pm$.

For $Q\in\Sigma$, define the formal sum (``generating function'')
\begin{equation}
\hat x(Q)=\sum_{n}x(n)\cdot \w^n(Q),\qquad
   Q\in\Sigma.                                \label{genf}
\end{equation}
From now on we always assume that while considering series of 1-forms or 2-forms
the index of summation runs over $\Z$ unless otherwise stated. Introduce the
operator-valued quadratic differential $T(Q)$ (the \emph{energy-momentum
tensor}) as follows:
\begin{equation}
T(Q):=\frac 12\sum_i\nord{\hat u_i(Q)\hat u^i(Q)} =
    \frac 12\sum_{n,m}\sum_i\nord{u_i(n)u^i(m)}\,\w^n(Q)\w^m(Q).
                              \label{emt}
\end{equation}
Here $\nord{\dots}$ denotes the normal ordering. Expand the quadratic
differential $T(Q)$ over the basis quadratic differentials:
\begin{equation}
T(Q)=\sum_kL_k\cdot\Omega^k(Q),
                                         \label{kazvect-3}
\end{equation}
where $L_k$ are operator-valued coefficients. Using duality relations we obtain
\begin{equation}
 L_k=\oint_{c_0} T(Q)e_k(Q)=\frac
     12\sum_{n,m}\sum_i \nord{u_i(n)u^i(m)}\,l_k^{nm},
                               \label{kazvect-4}
\end{equation}
where
\[
l_k^{nm}=\oint_{c_0} \w^n(Q)\w^m(Q)e_k(Q) .
\]
Notice that for a fixed $k$ the pairs $(n,m)$ such that $l_k^{nm}\ne 0$ satisfy
the inequality of the form $C_2\le m+n\le C_1$ where $C_1$, $C_2$ are constants
which depend only on $k$ and $g$.  Thanks to this property and to the normal
ordering, the $L_k$'s are well defined on admissible representations.  For
instance, if $m,n>0$ or $m,n<-g$ then $k\le m+n\le k+g$.  For $g=0$ this implies
$l_k^{nm}=\delta_k^{m+n}$ which leads to the usual definition of the Sugawara
operator for $g=0$ (see \cite{KaRa} and references therein).

We will consider the class of normal orderings which satisfy the following
requirement:
\begin{equation}
\nord{x(n)y(m)} = x(n)y(m),\quad \text{if  $n\le 0<m$ or $n<-g\le m$.}
                               \label{kazvect-5}
\end{equation}
The normal ordering determines the cohomology class of the cocycle in
Theorem~3.1 below.  In all other respects the following doesn't depend on the
choice of the normal ordering.

% Theorem 3.1
\begin{theorem}[\cite{SSS}]
Let $\g$ be a finite-dimensional Lie algebra, either abelian or simple. Let
$2\k$ be the eigenvalue of its casimir in the adjoint representation and $\gh$
the corresponding affine Krichever--Novikov algebra. Let $V$ be an admissible
representation of $\gh$ in which the central element acts as
$\ce\cdot\mathit{id}$.  If $\ce+\k\ne 0$, then the normalized Sugawara operators
$L_k^*=-(\ce+\k)^{-1}L_k$ define a representation of the Lie algebra $\Lh$ which
has a ``geometrical" cocycle
\[
\chi(e,f)=\frac {\ce\cdot\dim\g}{(\ce+\k)}\cdot \oint_{c_0}
    \biggl( \frac 12(e'''f-ef''')-R\cdot(e'f-ef')\biggr)\, dz ,
\]
where $e,f\in\L$, $R$ is a projective connection which has poles only
at the points $P_\pm$.
\end{theorem}

For each $e=\sum\l_ke_k\in\L$ (a finite sum) introduce $T(e)=\sum\l_kL_k^*$.  By
Theorem~3.1 $e\mapsto T(e)$ is a representation of $\Lh$. We call it the
\emph{Sugawara representation}.

In Section 4\ below we will need the following statement:

% Lemma 3.3
\begin{lemma}[\cite{SSS}]\leavevmode

\upn{(1)}\hspace{\labelsep}%
$[L_k,x(r)]=-(\ce+\k)\,x({e_k}A_r)$.

\upn{(2)}\hspace{\labelsep}%
$[L_k,\hat x(Q)]=(\ce+\k)e_k\ldot \hat{x}(Q)$,

\vspace{\topsep}
\noindent
where $e\ldot \hat{x}(Q):=\sum_n x(n)(e.\w^n)(Q)$, $eA$ is the derivative of
a function $A$, and $e.\w$ is the Lie derivative of a $1$-form $\w$ in the
direction of a vector field $e$.
\end{lemma}

% 3(d)
\subsection{The action of the Sugawara operators on the vacuum vectors}

% Lemma 3.4
\begin{lemma}
The Sugawara operators of the elements of the subalgebra $\L^{(2)}_+$ annihilate
the vacuum vectors of highest weight representations of $\gh$.
\end{lemma}

\begin{proof}
By Subsection (a), $e_k\in\L^{(2)}_+$ if and only if $k\ge 1$. Let us consider
the Sugawara operators $L_k$ for $k\ge 1$. We want to show that if $l_k^{mn}\ne
0$ then either $m\ge 1$ or $n\ge 1$. Provided it is true and taking into account
the normal ordering, one finds an operator of representation of the subalgebra
$\gh_+$ in the second position of the term $\nord{u(m)u(n)}$\,. Hence this term
annihilates the vacuum.

Consider the relation \eqref{kazvect-4} for $l_k^{mn}$. By \eqref{diff} and
\eqref{vect} $\ord_{P_+}\omega^m=-m-1$, $\ord_{P_+}\omega^n=-n-1$,
$\ord_{P_+}e_k=k+1$. Thus, $\ord_{P_+}\omega^m\omega^ne_k=-m-n+k-1$.  To have a
nonzero residue at the point $P_+$ the 1-form $\omega^m\omega^ne_k$ must have a
pole there, at least. Therefore $-m-n+k-1\le -1$. This implies $m+n\ge k\ge 1$
and hence either $m>0$ or $n>0$.
\end{proof}

% 4
\section{Casimirs, semi-casimirs, and moduli spaces}

In this section we give description of the second order casimirs for the Lie
algebra~$\gh$. Let $C_2$ denote the space of these casimirs. We also introduce
the semi-casimirs and establish their connection with the moduli space
$\M_{g,2}^{(p)}$ and coinvariants.

% 4(a)
\subsection{The second order casimirs}

For any affine Kac--Moody algebra there is only one second order
casimir\,---\,the sum of certain element of the Virasoro algebra and its
Sugawara operator \cite{KaRa}.  The main property of this operator is as
follows: it commutes with all the operators of the representation of the affine
algebra under consideration.

Based on this idea we will construct the second order casimirs for $\gh$ as
operators of the form $\Delta_e:=\eh -T(e)$, where $e\in\L$, $\eh$ is the
operator of a representation of the vector field $e$, and $T(e)$ is its Sugawara
operator. The words ``the second order" mean that the operators under
consideration depend quadratically on the operators of representation of the Lie
algebra $\gh$. Observe that even if we deal with the representation of
$\widehat{\mathfrak{gl}}_{g,2}$, the Sugawara representation $T$ corresponds to
the restriction of the latter onto the subalgebra
$\widehat{\mathfrak{sl}}_{g,2}$.

On the Riemann surface $\Sigma$, consider the Lie algebra $\D_\g$ of
differential operators of the form $e+xA$, $e\in\L$, $x\in\g$,
$A\in{\A}(\Sigma,P_\pm)$ (i.\,e., the algebra of Krichever--Novikov differential
operators of order less than or equal to~1).  As a linear space
$\D_\g=\L\oplus\g\otimes\A(\Sigma ,P_\pm)$. In particular, for
$\g={\mathfrak{gl}}(1)$ one has $\D_\g=\D$. The commutation relations between
vector fields and currents in $\D_\g$ are well known:
\begin{equation}
[e,x\otimes A]=x\otimes (eA).   \label{comm}
\end{equation}
Below, we will consider projective representations of $\D_\g$ (projective
$\D_\g$-modules). Such a representation is defined as a representation of
\emph{some} central extension $\widehat{\D_\g}$ of $\D_\g$. Assuming the action
of the central element to be the identity operator, we call the cocycle of this
central extension the \emph{cocycle of the projective ${\D_\g}$-module}
(\emph{representation}).  As it is stated in Section~2(d), the cocycle of the
fermion representation, while being restricted onto $\g\otimes\A(\Sigma,P_\pm)$,
gives the multiple of the cocycle \eqref{cocycle}. We call a projective
$\D_\g$-module possessing this property \emph{normalized}, and its cocycle as
well.  We call a projective representation of $\D_\g$ \emph{admissible} if its
restrictions to $\g\otimes\A(\Sigma ,P_\pm)$ and $\L$ are admissible.

% Lemma 4.1
\begin{lemma}
For a semisimple $\g$ and an admissible normalized projective $\D_\g$-module $V$
we have

$1^\circ$. $[\eh,x(A)]=x(eA)$ for any $A\in\A$, $e\in\L$.

$2^\circ$. $[\Delta_e, x(A)]=0$ for any $A\in\A.$
\end{lemma}

\begin{proof}
$1^\circ$. From \eqref{comm} and the definition of $\widehat{\D_\g}$-module, it
follows that
\begin{equation}
[\eh,x(A)]=x(eA)+\ga(e,xA)\circ id,                 \label{comm2}
\end{equation}
where $\ga$ is a cocycle on $\D_\g$.

For a semisimple $\g$ each cocycle on $\D_\g$ satisfies
\begin{equation}
\ga(e,xA)=0\qquad \text{for any}\ e\in\L,\ A\in\A,\ x\in\g.
                                                      \label{cocyc}
\end{equation}
Indeed, by definition of a cocycle
\[
\ga(e,[xA,yB])+\ga(yB,[e,xA])+\ga(xA,[yB,e])=0
\]
for any $e\in\L$, $A,B\in\A$, $x,y\in\g$.  The latter is equivalent to
\begin{equation}
\ga(e,[x,y]AB)+\ga(yB,x(eA))+\ga(xA,-y(eB))=0. \label{cocyc1}
\end{equation}
Take $B\equiv 1$ in the latter equality. Then $\ga(yB,x(eA))=\ga(y,x(eA))=0$,
since by the assumption of the Lemma $\ga$ is a multiple of \eqref{cocycle} and
the latter vanishes if one of the arguments is a constant.  Further on,
$eB\equiv 0$, therefore $\ga(xA,-y(eB))=0$.  The first summand in \eqref{cocyc1}
is equal to $\ga(e,[x,y]A)$.  Hence $\ga(e,[x,y]A)=0$ for all $e\in\L$,
$A\in\A$, $x,y\in\g$. If $\g$ is semisimple it coincides with its derived
algebra. Hence \eqref{cocyc} follows from the latter equality.

$2^\circ$ follows immediately from Lemma 3.3(1) and from $1^\circ$.
\end{proof}

Lemma 4.1.$2^\circ$ was proved in \cite{SSS} where our part $1^\circ$ above
appeared as a hypothesis.

% Lemma 4.2
\begin{lemma}[M.\,Schlichenmaier\footnotemark]
\footnotetext{Private communication.}  For a commutative $\g$ and an arbitrary
admissible representation $V$ of the corresponding \upn(Heisenberg type\upn)
affine algebra, we have
\[
[\Delta_e, x(A)]=\ga(e,A)\cdot id,\qquad\text{for any}\ e\in\L,\ A\in\A,
\]
where $\ga$ is a cocycle on $\D_\g$.
\end{lemma}

\begin{proof}
The relation \eqref{comm2} is always true (but this time $\ga$ can be a
nontrivial cocycle).  Comparing \eqref{comm2} with Lemma 3.3(1) proves the
claim.
\end{proof}

Now, let us consider the typical case of a reductive algebra.

%  Lemma 4.3
\begin{lemma}
For $\g={\mathfrak{gl}}(l)$ and $V$ and $e$ from Lemma~\upn{4.1}, the following
commutation relation holds\upn:
\begin{equation}
[\eh,x(A)]=x(eA)+\l(x)\ga(e,A)\circ \mathit{id},          \label{comm1}
\end{equation}
where $\ga$ is a cocycle on $\D$ and $\l(x)=l^{-1}\tr x$.
\end{lemma}

\begin{proof}
Again the relation \eqref{comm2} is true. An arbitrary $x\in\g$ can be
represented as follows: $x=x_0+\l(x)1_l$ where $x_0\in{\mathfrak {sl}}(l)$,
$\l\in\g^*$, $1_l$ stays for the identity matrix of rank $l$.  By Lemma~4.2,
$\ga(e,x_0A)=0$. Hence $\ga(e,xA)=\l(x)\ga(e,1_lA)$. Obviously the
correspondence $e+A\mapsto e+1_lA$ is a Lie algebra homomorphism
$\D\rightarrow\D_\g$. Hence $\ga(e,1_lA)$ defines a cocycle on $\D$.
\end{proof}

% Lemma 4.4
\begin{lemma}
Let $\g=\mathfrak{gl}(l)$ and let a representation of the corresponding algebra
$\gh$ be admissible.  Then for arbitrary $e\in\L$, $x\in\g$, $A\in\A$ we have
\[
[\Delta_e,xA]=\l(x)\ga(e,A)\circ \mathit{id},
\]
where $\ga$ is a cocycle on $\D$ and $\l(x)=l^{-1}\tr x$.
\end{lemma}

\begin{proof}
By Lemma 3 .3 one has $[T(e),x(A)]=x(eA)$.  Now the lemma follows from
comparison of the latter relation with Lemma 4.3.
\end{proof}

% Definition 4.1
\begin{definition}
The operators $\Delta_e$ which commute with all the operators of representations
of $\gh$ and $\A$ in $V$ are called (second order) \emph{casimirs} of the Lie
algebra $\gh$ in the representation $V$.
\end{definition}

% Remark 4.1
\begin{remark}
The requirement of commutativity of $\Delta_e$'s with $\A$ has a different
meaning for $\widehat{\mathfrak {sl}}_{g,2}$ and for
$\widehat{\mathfrak{gl}}_{g,2}$. Consider fermion representations as a typical
example. The space of monomials of a given charge is generically irreducible
with respect to $\widehat{\mathfrak{gl}}_{g,2}$ (see \cite{KaRa} for the
genus~$0$ case) but reducible with respect to $\widehat{\mathfrak{sl}}_{g,2}$.
Indeed, the elements of the form $d(\l)A\in\gh$ (here
$d(\l)=\mathit{diag}(\l,\dots,\l)$, where $\l\in\mathbb{C}$) belong to the
center of~$\gh$. They commute with all the operators of
$\widehat{\mathfrak{sl}}_{g,2}$ and hence make reducible its fermion
representation.  In this case the requirement of commutativity of $\Delta_e$'s
with $\A$ means that $\Delta_e$'s are well defined on the
$\widehat{\mathfrak{sl}}_{g,2}$-subrepresentations of the fermion
representation.
\end{remark}

The following lemma is an easy corollary of Definition 4.1 and Lemma 4.4.

% Lemma 4.5
\begin{lemma}
$\Delta_e$ is a casimir for $\gh$ \upn(in a certain representation\upn) if and only if
$\ga(e,A)=0$ for all $A\in\A$, where $\ga$ is the cocycle on $\D$ which
corresponds to the representation in question.
\end{lemma}

%\pdffix
% Remark 4.2
\begin{remark}
The following statement shows that being a casimir is a very restrictive
condition: \emph{given a cocycle $\ga$ on $\D$ the vector fields which satisfy
the condition of Lemma~4.5 form a Lie subalgebra in $\L$.}  To prove this
suppose $e_1,e_2\in\L$ to be such vector fields that $\ga(e_1,A)=\ga(e_2,A)=0$
for all $A\in\A$. By definition of a cocycle, $\ga([e_1,e_2],A)+\ga(e_2,[e_1,A])-
\ga(e_1,[e_2,A])=0$ for any two such vector fields and arbitrary $A$. Two terms
in the latter relation vanish because $[e_1,A]$, $[e_2,A]$ are functions, hence
$\ga([e_1,e_2],A)=0$.
\end{remark}

Denote the cocycle of the representation $V$ by $\ga_{V}$.  According to the
proof of Lemma~4.2, $\ga_{V}$ defines the cocycle $\ga_{V}(e,1_lA)$ on~$\D$.
We will keep the notation $\ga_{V}$ for this cocycle as well.

% 4(b)
\subsection{Casimirs of fermion representations}

Let us show that the fermion representations satisfy all the requirements of the
last subsection. The main point is that each fermion representation is an
admissible projective $\D$-module.

It was shown in Sections 2, 3 that a fermion representation is a projective
module over $\A$, $\g\otimes\A$ and $\L$.  The similar arguments show that it is
also a projective module over both $\D$ and $\D_\g$. Again let the fermion
representation be given by a holomorphic bundle $F$ with a meromorphic
(therefore flat) connection $\nabla$ which has logarithmic singularities at the
points $P_\pm$ (see Section~3(b)), and an irreducible representation $\tau$
of~$\g$.  By flatness, $[\nabla_e,\nabla_f]-\nabla_{[e,f]}=0$ for any
$e,f\in\L$.  Hence, $\nabla_{[e,f]}= [\nabla_e,\nabla_f]$ and $\nabla$ defines a
representation of $\L$ in the space of holomorphic sections of $F$.  By
definition of a connection, for each holomorphic section $s$ of the bundle, each
$e\in\L$ and each $A\in\A$ we have $\nabla_e(As) = (eA)s+A\nabla_es$ where $eA$
is defined by the relation \eqref{comm}.  In other words, $[\nabla_e,A]=eA$,
i.\,e., the mapping $e+A\rightarrow\nabla_e+A$ gives rise to a representation
of~$\D$.  Being conjugated by the matrix $\Psi$ (see Sections 2(b), 3(b)), this
representation can be transferred to the space $F_{KN}$ of Krichever--Novikov
vector-valued functions.  We will write down the action of the element $e\in\L$
on $\psi\in F_{KN}$ as $e\psi$. The space $V_F$, as introduced in Sections 2(c), 2(d), is
actually spanned by semi-infinite monomials over
$F_{KN}^\tau=F_{KN}\otimes V_\tau$ with the following $\D_\g$-action:
\[
(xA)(\psi\otimes v)=A\psi\otimes xv \quad\text{and}\quad
     e(\psi\otimes v)=e\psi\otimes v,
\]
where $V_\tau$ is the space of the representation $\tau$, $x\in\g$, $A\in\A$,
$\psi\in F_{KN}$ and $v\in V_\tau$.  Applying the procedure described above
(Sections~2,~3) of lifting representations from the space $F_{KN}^\tau$ to the
space $V_F$, we obtain projective representations of the algebras $\D$, $\D_\g$.
It is easy to show that these projective representations are admissible. So in
the case of fermion representations we are in the set-up of the previous
subsection.  In particular, we can study cocycles on $\D$ and $\D_\g$
arising from the fermion representations.

Let us consider the cocycles on $\D$ in fermion representations in more detail.
Take $\vac=\psi_M\wedge\psi_{M+1}\wedge\dots$ for a vacuum.  Let
\begin{equation}
\nabla_{e_k}=z^k \biggl(\partial +\omega_{-1}\frac{dz}{z}+O(1)\,dz\biggr) \label{nabloc}
\end{equation}
be the local behavior of $\nabla$ at $P_+$. For any $N\in\mathbb{Z}$, let $n(N)$
and $j(N)$ denote the first two components of the triple $(n,j,i)$ such that
$N=N(n,j,i)$ (Section~2(b)). Let $\omega_{j}$ be the $j$-th diagonal element of
$\omega_{-1}$.  Our main observation is the following

% Lemma 4.6
\begin{lemma}
For a fermion representation of $\D_\g$, assume that $M<0$, $\nabla$ satisfies
\eqref{nabloc} and $\sum_{N=M}^{-1}\omega_{j(N)}\notin\mathbb{Z}$.  The
cocycle $\ga$ of such a representation posesses the following properties\upn:
\begin{enumerate}
\item[$1^\circ$.] $\ga(A_{-k},e_k)\ne 0$, for any $k\in\mathbb{Z}$, $\k\ne 0$\upn;
\item[$2^\circ$.] $\ga(A_{-k},e_m)=0$, for any $m>k$.
\item[$3^\circ$.] $\ga(A_0,e_m)=0$, for any $m\in\mathbb{Z}$.
\end{enumerate}
\end{lemma}

\begin{proof}
$\ga(A_{-k}, e_m)=A_{-k}\circ \nabla_{e_m} -\nabla_{e_m}\circ A_{-k}-
e_mA_{-k}$.  This expression can be evaluated on the vacuum vector of the
representation.

$1^\circ$. For $m=k$ one of the summands $A_{-k}\circ \nabla_{e_k}$ or
$\nabla_{e_k}\circ A_{-k}$ annihilates the vacuum.  For example, for $k>0$ the
first summand does. Under assumptions of the lemma the $-\nabla_{e_k}\circ
A_{-k}-e_kA_{-k}$ doesn't anihilate the vacuum vector, as follows from direct
calculations with power series, for example, at $P_+$.

Locally, we have $A_{-k}=z^{-k}(1+O(z))$,
$e_k=z^{k+1}(1+O(z))\frac{\partial}{\partial z}$.  Let us act by these objects
on a Krichever--Novikov vector-valued function $\psi_N$ of degree $n$. While
considering the orders at the points $P_\pm$ one can forget for simplicity about
the conjugation of the vector field by the matrix $\Psi$. By \eqref{nabloc} we
have
\[
\nabla_{e_k}\circ A_{-k}\psi_N=
    \biggl(z^{k+1}\frac{\partial}{\partial z}+z^k\omega_{-1}\biggr)(z^{-k}z^n(1+O(z)) =
    (n-k+\omega_{-1})z^n(1+O(z)).
\]
For the other term we have $e_kA_{-k}= (-k)(1+O(z))$. Hence,
\[
(-\nabla_{e_k}\circ A_{-k}-e_kA_{-k})\psi_N=(2k-n+\omega_j)\psi_N+\dotsb,
\]
where dots denote the terms of higher degree, $n=n(N)$, $j=j(N)$.  By
regularization
\begin{equation}
(-\nabla_{e_k}\circ  A_{-k}-{e_k}A_{-k})\,\vac=
    \Biggl(\sum_{N=M}^{-1}(2k-n(N))+\omega_{j(N)}\Biggr)\vac .  \label{vac1}
\end{equation}
For $k<0$ the term $A_{-k}\circ \nabla_{e_k}$ enters the game instead of
$\nabla_{e_k}\circ A_{-k}$. This gives rise to
\begin{equation}
(A_{-k}\circ \nabla_{e_k}-{e_k}A_{-k})\,\vac=
    \Biggl(\sum_{N=M}^{-1}(k+n(N))+\omega_{j(N)}\Biggr)\vac.  \label{vac2}
\end{equation}
If $\sum_{N=M}^{-1}\omega_{j(N)}\notin\mathbb{Z}$ then both right hand sides of
\eqref{vac1}, \eqref{vac2} do not vanish for any $k\in\mathbb{Z}$.
Thus~$1^\circ$ is proven.

The proof of $2^\circ$, $3^\circ$ is similar.
\end{proof}

By Definition 4.1 and Lemma 4.4, $\Delta(e)$ is a casimir if and only if
\begin{equation}
\ga(A_k,e)=0, \quad \text{for any}\quad k\in\mathbb{Z}.  \label{sys}
\end{equation}
We will seek for the solutions to \eqref{sys} in the form
\begin{equation}
e=\sum_{m\ge m_0}a_me_m,                \label{form}
\end{equation}
where $m_0\in\mathbb{Z}$.

% Lemma 4.7
\begin{lemma}
For each fermion representation such that its cocycle $\ga$ satisfies conditions
of Lemma~\upn{4.6}, equation \eqref{sys} has one-dimensional space of solutions of the
form \eqref{form}.  For the generator of this space one has $m_0=0$.
\end{lemma}

\begin{proof}
For the vector fields of the form \eqref{form} the relations \eqref{sys} read as
\begin{equation}
\sum_{m\ge m_0}a_m\ga(A_{-k},e_m)=0, \quad \text{for all}\quad k\in\mathbb{Z},\ k\ne 0.
\label{sys1}
\end{equation}
This is an infinite system of linear equations for $a_m$'s. By Lemma 4.6 it has
a triangular matrix, hence, for each $k\ne 0$, $a_k$ can be expressed via the
$a_m$'s with $m<k$ from the $k$-th equation.  On the contrary, for $k=0$ we have
$\ga(A_0,e_m)=0$ for any $m\in\mathbb{Z}$, because $A_0=\mathit{const}$
(Lemma~4.6).  Thus $a_0$ is an independent constant. For $m<0$ one has $a_m=0$
because the latter is true for sufficiently large (negative) values of $m$. For
$m>0$ $a_m$'s can be expressed via $a_0$. This proves the claim.
\end{proof}

Bringing together Lemmas 4.5, 4.7 we obtain the following theorem.

% Theorem 4.1
\begin{theorem}
For each fermion representation such that its cocycle $\ga$ satisfies conditions
of Lemma~\upn{4.6} there exists exactly one \upn(up to a scalar factor\upn)
casimir.  The corresponding vector field has a simple zero at $P_+$.
\end{theorem}

Let us normalize a casimir by the condition $a_0=1$. Then the eigenvalue of this
casimir equals to the one of $\Delta(e_0)$, because by Lemmas 3.2, 3.4
$\Delta(\L^{(2)}_+)$ annihilates the vacuum.

For an arbitrary vector field of the form \eqref{form} we can claim that both
the action of the vector field and its Sugawara operator are well defined. To
see the first, let us recall that for each $v\in V$ and $n$ sufficiently large,
$e_nv=0$.  For Sugawara operators one has $L_e=\sum_{m,n}l^{m,n}_e\,\nord{u_mu_n}$,
where $l^{m,n}_e=\res_{P_+}(\omega^m\omega^n e)$. But for given $m,n$
$\res_{P_+}\omega^m\omega^ne_k\ne 0$ only for a finite number of $k$'s. Hence
the coefficients $l^{m,n}_e$ are also well defined.

% 4(c)
\subsection{Semi-casimirs and moduli spaces $\M_{g,2}^{(p)}$}

From Lemma 4.5 and the previous subsection one can conclude that the condition
for $\Delta_e$ to commute with all the elements of the algebra $\A$ imposes very
strong restrictions on~$e$.  Here we consider weaker conditions.

% Definition 4.2
\begin{definition}
We call an operator of the form $\Delta_e$ a semi-casimir if $[\Delta_e,A_k]=0$
for each $k<0$.
\end{definition}

Let $\At_-\subset\A$ be the subspace spanned by all $A_k$, $k<0$.  Observe that,
as a Lie algebra, $\At_-$ is a subalgebra of $\A$ (while, as an associative
algebra, it is not). For $\g=\mathfrak{gl}(l)$, $\At_-$ can be also considered
as a Lie subalgebra of $\g\otimes\A$ (but not of $\gh$) by interpreting elements
of $\A$ as scalar matrices. As a subalgebra of $\g\otimes\A$, $\At_-$ commutes
with $\mathfrak{sl}(l)\otimes\A_-$.  Hence $\g_r=\mathfrak{sl}(l)\otimes\A_-
\oplus\At_-$ is a Lie subalgebra of $\g\otimes\A$. We call $\g_r$ a
\emph{regular subalgebra} (note that a different choice of $\g_r$ is possible,
see \cite{Shec}).

Define \emph{coinvariants} of the regular subalgebra in a $\gh$-module
$V$ as the quotient space $V/\U(\g_r)V$, where $\U(\g_r)$ denotes the
subalgebra of the universal enveloping algebra of $\gh$ corresponding
to the subspace $\g_r$ (in particular, it does not contain the unit).

By Definition 4.2, $\Delta_e$ is well defined on the space of coinvariants
of~$\g_r$.  Thus, this kind of condition is important for conformal field
theory.  On the other hand, Remark~4.2 is not applicable to semi-casimirs,
because the corresponding vector fields do not form a subalgebra but only a
linear subspace.

For a vector field $e$ defining a semi-casimir one has
\begin{equation}
\ga(A_{-k},e)=0,\quad \text{for any $k\in\mathbb{Z}$, $k>0$,}  \label{con2s}
\end{equation}
instead of \eqref{sys}.  For semi-casimirs one has system of equations similar
to \eqref{sys1} but only for $k>0$.  Therefore, the coefficients $a_m$ with $m\le 0$
appear to be independent and all the others can be expressed via them. Let
$\Lt_-\subset\L$ be the subspace spanned by $\{e_k\colon k\le 0\}$.  Introduce
the map $\Gamma\colon\Lt_-\mapsto\L$ as follows: take $e\in\Lt_-$ and represent
it in the form \eqref{form}; then substitute the corresponding $a_m$ ($m\le 0$)
into \eqref{sys1} and calculate $a_m$, $m>0$. Denote by $\Gamma(e)$ the vector
field which corresponds to the full set of $a_m$'s. We claim the following.

% Lemma 4.8
\begin{lemma}
The space of semi-casimirs coincides with $\Delta(\Gamma(\Lt_-))$. It is spanned
by the elements $\Delta(\Gamma(e_k))$, where $k\le 0$.
\end{lemma}

As it was mentioned above, semi-casimirs are well defined on the space of
coinvariants of the subalgebra $\g_r$. For $e\in\Lt_-$, let
$\overline{\Delta}(e)$ be the operator induced by $\Delta(\Gamma(e))$ on
coinvariants.  The map $\overline{\Delta}$ is defined on $\Lt_-$ and by
Lemma~4.8 its image is the space $C^s_2$ of semi-casimirs considered as
operators on the space of coinvariants.

Our next step is to show that only a finite number of basis semi-casimirs are
nonzero on coinvariants and, for a proper natural $p$, to establish the
correspondence between the tangent space to $\M _{g,2}^{(p)}$ and the space of
semi-casimirs (considered on coinvariants).

% Lemma 4.9
\begin{lemma}
For a fermion representation $V$ there exists such $p\in\mathbb{Z}_+$ that
$\L^{(p)}_-\subseteq\ker\overline{\Delta}$.
\end{lemma}

\begin{proof}
For a fermion representation of certain charge, the space of coinvariants of
$\g_r$ is finite-dimensional. It is a known fact for affine Kac--Moody
algebras. It is easy to reduce the statement in the almost graded case in
question to the known case just by considering the associated graded objects.

Consider the decomposition
\[
V=\bigoplus\limits_{n\le 0}V_n
\]
defining the structure of an almost graded $\D_\g$-module on $V$.
By finite-di\-men\-sion\-al\-ity of coinvariants there exists such $s$ that
\[
\bigoplus\limits_{n<s}V_n\subseteq\U(\g_r)V.
\]
Since $V$ is an almost graded $\L$-module there exists such $\nu\in\mathbb{Z}_+$
that $\eh_kV_n\subseteq\bigoplus_{m\le k+n+\nu}V_m$. Obviously, $k+n+\nu<s$ for
any $n\le 0$ if $k<s-\nu$.  Hence, $\eh_kV\subseteq \U(\g_r)V$ for any
$k<s-\nu$.

Let us find such $k'\in\mathbb{Z}$ that $T(e_k)V\subseteq \U(\g_r)V$ for any
$k<k'$. Since $V$ is an almost graded $(\g\otimes\A)$-module there exists such
$\nu'\in\mathbb{Z}_+$ that $u(i)V_n\subseteq\bigoplus_{m\le i+n+\nu'}V_m$ for
any $u\in\g$, $i\in\mathbb{Z}$. Hence, $u^{\mu}(i)u_{\mu}(j)V_n\subseteq
\bigoplus_{m\le i+j+n+2\nu'}V_m$.  The term $\nord{u^{\mu}(i)u_{\mu}(j)}$ occurs in
the series \eqref{kazvect-4} for $T(e_k)$ if $l^{ij}_k\ne 0$. The latter holds
for
\begin{equation}
i+j\le k+g.         \label{non}
\end{equation}
To obtain this observe that, assuming $g>1$,  for the basis elements
$A_m\in\A$, $\omega_m$ (the 1-form), $e_m\in\L$, one has the following behaviour at
$P_-$: $A_m=O(z^{-m-g})$, $\omega_m=O(z^{m+g-1})dz$,
$e_m=O(z^{-m-3g+1})\frac{\partial}{\partial z}$.  These relations also hold for
$g=1$ if one sets $\omega_m=A_{1-m}dz$, $e_m=A_{m+1}\frac{\partial}{\partial
z}$, and also for $g=0$.  By \eqref{kazvect-4}
$l_k^{ij}=-\res_{P_-}e_k\omega^i\omega^j$. At the point $P_-$ one has
$e_k\omega^i\omega^j=O(z^{-k+i+j-g-1})$. If $l^{ij}_k\ne 0$ then $-k+i+j-g-1\le
-1$, which implies \eqref{non}. Therefore,
\[
T(e_k)V_n\subseteq\bigoplus\limits_{m\le n+k+g+2\nu'}V_m.
\]
Hence, we can take $k'=s-g-2\nu'$.

Obviously, for $p=\max(\nu-s,|k'|)$ and $e\in\L^{(p)}_-$ we have
$\Delta(e)V\subseteq\U(\g_r)V$, thus $e\in\ker\overline{\Delta}$.
\end{proof}

Let $\M _{g,2}^{(p)}$ be the moduli space of curves of genus $g$ with two marked
points $P_\pm$, fixed $1$-jet of local coordinate at $P_+$ and fixed $p$-jet of
local coordinate at $P_-$.  There is a canonical mapping $\theta:\L\mapsto
T_\Sigma\M _{g,2}^{(p)}$. This mapping goes back to \cite{Kon}.  The
cohomological and geometrical versions of this mapping are given in \cite{SSkz}
and \cite{GrOr}, respectively. Let $\tilde\theta$ denote the restriction of
$\theta$ onto the subspace $\Lt_-$.  Let $V$ be a fermion representation of
$\D_\g$ and $\ga_V$ be the cocycle of this representation.  Let also
$C_2^s=C_2^s(V)$ denote the second order semi-casimirs of $\gh$ in the
representation $V$. We assume semi-casimirs to be restricted onto coinvariants.

% Theorem 4.2
\begin{theorem}
Take $p$ as in Lemma \upn{4.9}.
\begin{enumerate}
\item[$1^\circ$.]  The mapping $\tilde\theta\colon \Lt_-\mapsto
T_\Sigma\M_{g,2}^{(p-1)}$ is surjective and\/ $\ker\tilde\theta=\L^{(p)}_-$.

\item[$2^\circ$.]  For such $V$ that $\ga_V$ satisfies the conditions of
Lemma~\upn{4.6}, the mapping $\overline{\Delta}\colon \Lt_-\mapsto C_2^s(V)$ is
surjective and $\L^{(p)}_-\subseteq\ker\overline{\Delta}$.

\item[$3^\circ$.] The mapping $\overline{\Delta}\circ\tilde\theta^{-1}\colon
T_\Sigma\M_{g,2}^{(p-1)} \mapsto C_2^s(V)$ is surjective.
\end{enumerate}
\end{theorem}

\begin{proof}
$1^\circ$ is an easy corollary of the known results just mentioned.  By
Proposition 4.4 of \cite{SSkz}, the mapping $\theta\colon\L\mapsto
T_\Sigma\M_{g,2}^{(p-1)}$ is surjective and $\ker\theta=\L^{(2)}_+\oplus
\L^{(p)}_-$.  The claim follows due to the decomposition
$\L=\L^{(2)}_+\oplus\Lt_-$.

$2^\circ$ follows from Lemmas 4.8, 4.9.  Part $3^\circ$ is an immediate
consequence of $1^\circ$ and~$2^\circ$.
\end{proof}

% 5.
\section{Casimirs in terms of ``geometrical" cocycles}

The goal of this section is to repeat some of the results of the previous
section in another set-up. A 2-cocycle on $\D$ is called local if there exists
$L\in\mathbb{Z}$ such that $\ga(A_i,A_j)=\ga(e_i,e_j)=\ga(A_i,e_j)=0$ for any
$i,j\in\mathbb{Z}$ such that $|i-j|>L$. The local cocycles on $\D$ have the
following description.

% Lemma 5.1
\begin{lemma}
Any local cocycle on $\D$
in local coordinates is of the form
\begin{multline*}
    \ga(f_1\partial +g_1,f_2\partial+g_2)
    =\res_{P_+}
    [a_1(f_1f_2'''-f_1'''f_2) + R(f_1f_2'-f_1'f_2) +
    a_2(f_1g_2''-f_2g_1'')\\
    +T(f_1g_2'-f_2g_1') + a_3g_1dg_2],
\end{multline*}
where $a_1,a_2,a_3\in\mathbb{C}$.
\end{lemma}

For $g=0$ this statement and its proof are contained in the proof of the
Proposition~2.1 form \cite{ArCoKaPr}.  For $g>1$ the author doesn't know any
published proof of Lemma~5.1 (see Introduction). Nevertheless the statement is
wellknown; it gives rise to another approach to (semi-)casimirs. Here, we
consider this approach.

% Lemma 5.2
\begin{lemma}
For any local cocycle $\ga$ on $\D$, the following relation in local coordinates
holds\upn:
\[
\ga(e,A)=\res_{P_+}(afA''+TfA'),
\]
where $e=f\partial$, $a\in\mathbb{C}$, the behavior of $T$ under an
arbitrary change of local coordinates can be described by any of the
following three equivalent relations\upn:
\begin{enumerate}
\item[$1^\circ$.]
         $a^{-1}T(u)=a^{-1}T(z)u_z^{-1}+u_{zz}u_z^{-2}$\upn;
\item[$2^\circ$.]
        $a^{-1}T(u)du=a^{-1}T(z)dz+d\ln u_z$\upn;
\item[$3^\circ$.]
        there exists a $v$ such that $T(z)=a\frac{\partial}{\partial z}\ln v(z)$,
\end{enumerate}
where $u,z$ are local parameters, $u_z=u'(z)$, $v\partial\in\L$ is the local
representation of a vector field.
\end{lemma}

\begin{proof}
The expression for $\ga$ can be obtained immediately by taking $f_1=f$, $g_1=0$,
$f_2=0$, $g_2=A$ in Lemma 5.1.

Now, let us find the transformation law that a pair $a,T$ must satisfy in order
to make the expression $\ga_{a,T}(e,A)= \res_{P_+}(afA''+TfA')$ independent of
the choice of local coordinates. This independence is obviously equivalent to
the requirement that $\Omega_{a,T}(A):=aA''+TA'$ be a quadratic differential for
each $A\in\A$. It is easy to show that the desired transformation law can be
written in the three equivalent forms above. For example let us briefly outline
the implication $2^\circ\,{\Rightarrow}\,3^\circ$. By integration of both parts of
the equality $2^\circ$ from any fixed point to the current point $P$ one obtains
$a^{-1}\int^PT(u)du=a^{-1}\int^PT(z)dz+\ln u_z$. By exponentiating, the latter
equality implies $\exp(a^{-1}\int^PT(u)du)=\exp(a^{-1}\int^PT(z)dz)u_z$.  Now it
is easy to notice that $v(z)=\exp(a^{-1}\int^PT(z)dz)$ transforms like a vector
field.
\end{proof}

% Remark 5.1
\begin{remark}
Part $3^\circ$ of Lemma~5.2 makes sure that the quantities $T$ with the required
transformation law do exist.  One can take an arbitrary vector field from $\L$
and construct $T$ as prescribed by the relation in part $3^\circ$.
\end{remark}

Using Lemma 5.2, Lemma 4.4 can be made more precise as follows.

% Lemma 5.3
\begin{lemma}
Let $\g=\mathfrak{gl}(l)$ and suppose the representation of the corresponding
algebra $\gh$ is admissible.  Then for arbitrary $e\in\L$, $x\in\g$, $A\in\A$
one has
\[
[\Delta_e,xA]=\l(x)\ga(e,A)\circ\mathit{id},
\]
where in local coordinates for $e=f\partial$ the following relation holds\upn:
$\ga(e,A)=\res_{P_+}(afA''+TfA')$ \upn($a\in\mathbb{C}$\upn) and $T$ satisfies
the conditions of Lemma~\upn{5.2}.
\end{lemma}

From Lemmas 5.6, 5.3 one can see that only the cocycles
of the form
\[
\ga_{a,{}_T}(f_1\partial +g_1,f_2\partial+g_2)=
\res_{P_+}[a(f_1g_2''-f_2g_1'')+T(f_1g_2'-f_2g_1')],\quad a\in\mathbb{C}^\times,
\]
are responsible for commutativity of the operators $\Delta_e$ and $x(A)$.

By Lemma~5.2 there exist $T_V$ and $a\in\mathbb{C}$ such that $\ga(e,A)=
\res_{P_+}(afA''+T_VfA')$ and $aA''+T_VA'$ is a quadratic differential for each
$A\in\A$. Let $\Omega_V$ denote the linear space of all quadratic differentials
of this form.  Introduce also the notation $\Omega^{(2)}$ for the space of all
meromorphic quadratic differentials and $\langle\,\cdot\,{,}\cdot\,\rangle$ for the
natural pairing between vector fields and quadratic differentials: $\langle e,
\Omega\rangle:=\res_{P_+}(e\Omega)$ ($e\in\L$, $\Omega\in\Omega^{(2)}$).

% Theorem 5.1
\begin{theorem}
$\Delta_e$ is a casimir of $\gh$ in the representation $V$ if and only if
$\langle e,\Omega\rangle=0$ for each $\Omega\in\Omega_V$.
\end{theorem}

\begin{proof}
By definition of $\langle\,\cdot\,{,}\cdot\,\rangle$, Lemma 5.3 reads as
\[
[\Delta_e,xA]=\l(x)\langle e,\Omega\rangle \circ\mathit{id},
\]
where $\Omega=aA''+T_VA'$.  It can be easily extracted from the proof of
Lemma~5.3 that $\l(1_l)=1$, hence,
\[
[\Delta_e,A]=\langle e,\Omega\rangle \circ\mathit{id}.
\]

These two relations mean that $\Delta_e$ commute with all the elements of the
form $A$ and $xA$ ($x\in\g$, $A\in\A$) iff $\langle e,\Omega\rangle=0$ for each
$\Omega=aA''+T_VA'$. If $A$ runs over $\A$ then $\Omega$ runs over
$\Omega_V$. This proves the theorem.
\end{proof}

Let $\Omega_V^\bot :=\{e\in\L\colon \langle e,\Omega\rangle=0, \forall
\Omega\in\Omega_V \}$.

% Corollary 5.1
\begin{corollary}
$C_2\cong \Omega_V^\bot/(\Omega_V^\bot\cap \L^{(2)}_+)$.
\end{corollary}

This follows from Theorem 5.1 and Lemmas 3.2, 3.3.

% 5(b)
\subsection{Description of casimirs}

In a local coordinate $z$ on $\Sigma$ a vector field $e$ can be represented in
the form $e=E\partial$ where $E=E(z)$ is a local function,
$\partial=\frac{\partial}{\partial z}$.  In a generic situation $T$ can be
chosen to have simple poles at the points $P_\pm$ (Lemma 5.2.$3^\circ$).
Example~5.1 below shows that this is the case also for those a$T$'s which come
from fermion representations. Assume that a representation $V$ is fixed and
$T=T_V$.

By Theorem 5.1 the condition for $\Delta_e$ to be a casimir reads
\begin{equation}
\oint_{c_0} E(aA''+TA')dz =0,\quad \text{for any}\ A\in\A.  \label{cond1}
\end{equation}
Integrating ``by parts'' and taking into account the arbitrariness of $A$ one
obtains the differential equation
\begin{equation}
aE''-(ET)' =0.                             \label{cond3}
\end{equation}
We will seek the solutions to \eqref{cond3} of the form
\begin{equation}
E=\sum_{n\ge N}a_nE_n,                \label{formd}
\end{equation}
where $e_n=E_n\frac{\partial}{\partial z}$ near the point $P_+$.

% Lemma 5.4
\begin{lemma}
For a generic $T$ such that $T(z)=O(z^{-1})$ at the point $P_+$ the equation
\eqref{cond3} has the one-dimensional space of solutions of the form
\eqref{formd}.
\end{lemma}

\begin{proof}
The proof is nothing but expanding equation \eqref{cond3} in power series in
the neighbourhood of the point $P_+$.  Suppose for simplicity that
$\ord_{P_+}e=-1$.  By assumptions we have $E=\ep_{-1}z^{-1}+\ep_0+{\ep_1}z+
\dotsb$, $T=\tau_{-1}z^{-1}+\tau_0+{\tau_1}z+\dotsb$. Hence, for the power series
the relation \eqref{cond3} reads as
\begin{equation}
\begin{aligned}
     -2\ep_{-1}(1+\tau_{-1})&=0,\\
     \ep_0\tau_{-1}+\ep_{-1}\tau_0&=0,\\
     \ep_1\tau_0-(2-\tau_{-1})\ep_2&=0,\\
      \omit\span\omit\leaders\hbox{\,.\,}\hfil
\end{aligned}
\label{coef}
\end{equation}
For a generic $T$ one has $1+\tau_{-1}\ne 0$, hence the first relation implies
$\ep_{-1}=0$. Similarly $\ep_0=0$. From the third relation we obtain
$\ep_2=\frac{\tau_0}{2-\tau_{-1}}\ep_1$. So we have exactly one independent
constant $\ep_1$. All the remaining relations express $\ep_k$'s ($k>\nobreak1$) via
$\ep$'s with smaller numbers. This remains true if $\ord_{P_+}e<-1$.
\end{proof}

The following theorem is an immediate consequence of Lemma 5.4.

% Theorem 5.3
\begin{theorem}
For a generic $T$ such that $T(z)=O(z^{-1})$ at the point $P_+$, there is only
one \upn(up to a scalar factor\upn) second order casimir. This casimir corresponds
to the vector field which behaves at the point $P_+$ as follows\upn:
$e(z)=z(1+O(z))\frac{\partial}{\partial z}$.
\end{theorem}

% 5(c)
\subsection{Description of semi-casimirs}

By Theorem 5.1 and in analogy with \eqref{con2s}, for a vector
field which defines a semi-casimir one has
\begin{equation}
\res_{P_+} (aE''-(ET)')A\,dz=0,\quad \text{for any}\ A\in\At_-. % \label{con2s}
\end{equation}
Let us suppose $a=1$ for simplicity and take $E''-(ET)'=\sum a_iz^i$ at the
point $P_+$ (the sum is finite from the left). We deal only with local
expansions of the objects at the point $P_+$ here. For the element
$A_{-1}=\a_{-1}z^{-1}+ \a_0+\dotsb$ of the highest degree in $\At_-$, the
relation \eqref{con2s} gives $\b_{-1}a_0+\b_0a_{-1}$ which enables one to
express $a_0$ via $a_{-1}$. Similarly for the next basis element
$A_{-2}\in\At_-$, \eqref{con2s} expresses $a_1$ via $a$'s with the smaller
numbers. Thus the coefficients $a_i$, $i<0$ are independent and all the others
can be expressed via them.

We have the following relations for the semi-casimirs instead \eqref{coef}:
\begin{equation}
\begin{aligned}
     \omit\span\omit\phantom{$\ep_{-1}\tau_2+\ep_0\tau_1+\ep_1$}
    \leaders\hbox{\,.\,}\hfil\\
     -2\ep_{-1}(1+\tau_{-1})&=a_{-3},\\
     \ep_0\tau_{-1}+\ep_{-1}\tau_0&=a_{-2}, \\
     \ep_{-1}\tau_2+\ep_0\tau_1+\ep_1\tau_0+\ep_2\tau_{-1}-2\ep_2&=a_0,\\
      \omit\span\omit\leaders\hbox{\,.\,}\hfil
\end{aligned}
\label{coefs}
\end{equation}
The first two relations express $\ep_{-1}$, $\ep_0$ via independent constants
$a_{-3}$, $a_{-2}$. The third relation expresses $\ep_2$ via $\ep_1$.  The other
relations express $\ep_i$'s via the $\ep$'s with smaller numbers ($i\ge 2$). So
the $\ep_i$'s are independent if and only if $i\le 1$. Thus we have obtained
Lemma~4.8 again.

The following example shows that highest weight representations which satisfy
the conditions of Theorem 5.3 actually exist.

% Example 5.1
\begin{example}
Consider the fermion representation which corresponds to a generic 2-dimensional
bundle on an elliptic curve.  We are going to show that the cocycle of this
representation corresponds to the affine connection $T$ which has the pole of
order~1 at the point $P_+$.

Let $\ga$ stand for this cocycle. First of all let us show that \emph{if
$\ord_{P_+}T<-1$ then there exists $m\ge 0$ such that $\ga(A_1, e_m)\ne 0$}.
Suppose $\ord_{P_+}T=-m-2$, $m\ge 0$. Since $\ord_{P_+} e_m=m+1$ and both
$A_1''$ and $A_1'$ are holomorphic ($A_1'(P_+)\ne 0$ for a generic situation),
one has $(A_1''+TA_1')e_m=O(z^{-1})dz$. Hence, for a generic situation
$\res(A_1''+TA_1')e_m\ne 0$.

But it is easy to show for a highest weight representation that $\ga(A_1,e_m)=0$
for each $m\ge 0$. One has $[A_1,e_m]=-e_mA_1+\ga(A_1,e_m)$. For an elliptic
curve one can take $e_m=A_{m+1}\frac{\partial}{\partial z}$. Hence,
$\ga(A_1,e_m)=-[A_1,e_m]-e_mA_1=e_m\circ A_1-A_1\circ e_m- A_{m+1}
\frac{\partial A_1}{\partial z}$. Apply both parts of the latter equality to the
vacuum. Then one has $A_1\,|0\rangle =0$, $e_m\,|0\rangle =\l\,|0\rangle$ where
$\l\in\mathbb{C}$. Evidently $A_{m+1}\frac{\partial A_1}{\partial z}$ is an
element of the subalgebra $\A_+$. Hence, $\ga(A_1, e_m)= 0$.  What we have
obtained is
\[
\ord_{P_+}T\ge -1.
\]

Therefore we have to decide which of the two possibilities takes place: $T$ is
regular at the point $P_+$ or $T$ has a pole of order~$1$.  It is easy to check
that in the first case $\res_{P_+}(A_1''+TA_1')e_{-1}=0$ while in the second
case this residue generically doesn't vanish. Hence, loking at the value of
$\ga(A_1, e_{-1})$ we can understand which one of the two above possibilities
takes place.  By Lemma~4.6 generically $\ga(A_1, e_{-1})\ne 0$ hence,
$\ord_{P_+}T=-1$.
\end{example}

% 6.
\section{Semi-casimirs and quantization of the second order
              Hitchin integrals}

In this section we non-formally show how the semi-casimirs 
appear in course of operator quantization of the second order Hitchin 
integrals.

The second order Hitchin integrals $\chi_i$ are defined by the
expansion $tr\,\phi^2=\sum \chi_i\Omega^i$ where $\phi$ is the Higgs
field and $\{\Omega^i\}$ is a base of the cotangent
space to ${\mathcal M}_{g,2}^{(p-1)}$ realized as a certain
space of quadratic differentials.

In course of operator quantization the Higgs field $\phi$ should be
replaced by the current $I$. $I$ is an arbitrary Krichever-Novikov
1-form (on the Riemann surface) which has values in the
representation operators of $\gh$. Thus $I=\sum u_k\omega^k$ where
$\omega^k$ are the basis Krichever-Novikov 1-forms, $u_k$ are the
operator-valued coefficients ($k\in{\Bbb Z}$). The $\phi^2$ should be
replaced by $\nord{I^2}$.  The trace $tr\,\nord{I^2}$ (where $tr$ is
the "finite-dimensional" trace which means that it is linear over the
function algebra $\A$) is just the energy-momentum tensor $T$ which
was introduced in Section~3.  The expansion $T=\sum L_i\Omega^i$ (cf.
\eqref{kazvect-3}) is the quantum analog of the above expansion
$tr\,\phi^2=\sum \chi_i\Omega^i$. We will consider the normalized
form $-T(e_i)$ of an operator $L_i$ (see Section~3). What is usually
being done for compensating a
normal ordering is adding certain cartanian elements to the normal
ordered quantity.  For example, for Kac-Moody algebras the vector
field $z\frac{\partial}{\partial z}$ is being added, i.e. the
casimir equals to $z\frac{\partial}{\partial
z}-T(z\frac{\partial}{\partial z})$.  Observe that making use of
this idea exactly results in semi-casimirs $\Delta_i=e_i-T(e_i)$.
Since there is only finite number of Hitchin integrals and the
infinite set of semi-casimirs $\Delta_i$ we must formulate some
selection rule for the latters.  We propose to consider only those
$\Delta_i$ which induce nontrivial operators on conformal blocks.
Then by Theorem~4.2 we obtain the natural mapping of the $\chi_i$'s
to the $\Delta_i$'s.

\bibliographystyle{amsalpha}
%\bibliography{sheinman}

\providecommand{\MR}{\relax\ifhmode\unskip\space\fi MR }
% \MRhref is called by the amsart/book/proc definition of \MR.
\providecommand{\MRhref}[2]{%
  \href{http://www.ams.org/mathscinet-getitem?mr=#1}{#2}
}
\providecommand{\href}[2]{#2}

\end{document}